\documentclass[a4paper]{amsart}
\usepackage{amsmath,amsthm,amssymb,latexsym,epic,bbm,comment,color}
\usepackage{graphicx,enumerate,stmaryrd}
\usepackage[all,2cell]{xy}
\xyoption{2cell}

\usepackage{tikz-cd}

\renewcommand{\epsilon}{\varepsilon}

\newtheorem{theorem}{Theorem}

\newtheorem{definition}[theorem]{Definition}
\newtheorem{corollary}[theorem]{Corollary}
\newtheorem{proposition}[theorem]{Proposition}

\newtheorem{example}[theorem]{Example}

\newtheorem{remark}[theorem]{Remark}

\begin{document}
	\title[A Koszul-Type Splitting Theorem for Complex Supermanifolds]{A Koszul-Type Splitting Theorem for Complex Supermanifolds}

	\author{Luiz F.A. Campos, Alexey Sharapov and \\
    Elizaveta Vishnyakova}

\begin{abstract}

We study the splitting problem for complex supermanifolds using formal exponential maps 
and Ehresmann connections. We show that every formal exponential map induces a splitting 
of the underlying supermanifold. Furthermore, we establish that every Ehresmann connection 
determines a linearization of the infinite jet space, extending the corresponding construction 
for ordinary manifolds to supermanifolds.

For Ehresmann connections preserving the zero section, we construct a formal exponential map 
in two ways: via the linearization of the infinite jet space and directly from the connection 
operator. We prove that both constructions yield the same map and induce the same splitting, 
thereby generalizing Koszul's splitting theorem from linear connections to zero-section-preserving 
Ehresmann connections. Finally, we provide examples showing that, without this assumption, 
the two constructions may yield different splittings, or one or both of them may fail to induce 
a splitting.

\end{abstract}

	\maketitle

\section{Introduction}

The splitting problem is one of the central questions in complex supergeometry. A complex supermanifold $\mathcal{M}$ is called \emph{split} if it is isomorphic to its associated graded supermanifold $\operatorname{gr}\mathcal{M}$. While every smooth supermanifold is split by Batchelor-Gawedzki theorem, this is famously not the case in the complex analytic category, where non-split supermanifolds are generic. The splitting problem and its cohomological obstructions have consequently attracted considerable attention; see, for example, \cite{Correa2026-gm,Vi} and references therein. 

It is therefore natural to investigate geometric structures that ensure the existence of a splitting. A foundational result in this direction is Koszul's theorem \cite{Kos}, which guarantees a splitting in the presence of a holomorphic linear connection. The main purpose of the present paper is to generalize Koszul's theorem to a much broader class of connections, utilizing the geometry of formal exponential maps and infinite jet spaces.

We begin by investigating Ehresmann connections on complex supermanifolds (specifically, on their tangent supermanifolds) and their relation to infinite jet spaces. We prove that every Ehresmann connection determines a canonical linearization of the infinite jet space, compatible with the natural projections between jet spaces of different orders. This construction extends to the supergeometric setting the classical linearization known in ordinary differential geometry. Notably, no assumption regarding the preservation of the zero section is required for this result.

Next, we introduce formal exponential maps on complex supermanifolds,
extending the definition given in \cite{CF} to the supergeometric setting. 
Formal exponential maps were also studied by Emmrich and Weinstein~\cite{Emmrich1993TheDG} in the context of Fedosov's deformation quantization on ordinary manifolds, and this construction was subsequently extended to supermanifolds in~\cite{LS}. In the present work, we show that every formal exponential map induces a splitting of the underlying supermanifold. The proof uses a natural quotient of the sheaf of formal functions on the tangent supermanifold, which can be identified with the structure sheaf of the associated graded supermanifold. This quotient is closely related to the construction considered in \cite{Vi}.

The core of the paper is devoted to constructing formal exponential maps from Ehresmann connections. We develop two distinct approaches: the first is based on the linearization of infinite jet spaces, while the second employs the formal exponential of the differential operator associated with the connection. Our main theorem asserts that if the Ehresmann connection preserves the zero section, both constructions yield well-defined, normalized formal exponential maps that strictly coincide. Consequently, they induce the same splitting of the supermanifold. Since linear connections are particular instances of zero-section-preserving Ehresmann connections, this result directly generalizes Koszul's splitting theorem~\cite{Kos}.

Finally, we address Ehresmann connections that do not preserve the zero section. In this setting, the exponential constructions are not guaranteed to be well-defined as formal power series. Nevertheless, in certain cases, their coefficients can be summed analytically, yielding well-defined geometric morphisms. We exhibit explicit examples demonstrating that, in the absence of this condition, the two constructions may produce distinct splittings, only one of them may induce a splitting, or both may fail to yield a splitting altogether.

\section*{Acknowledgments:} 
E.V. was partially  supported 
by FAPEMIG grant RED-00133-21 - Rede Mineira de Matem\'{a}tica, CNPq grant 402320/2023-9. The results of Section~4  
were obtained under the exclusive support of the Ministry of Science and Higher Education of the Russian Federation (project FSWM-2025-0007). 
The authors would like to thank Alejandro Cabrera, Yuri Conrado, and 
Luca Vitagliano for helpful discussions related to this work.

\section{Preliminaries}

\subsection{Jet coordinates on a supermanifold}\label{sec:jet-coordinates} 

We work in the category of smooth or complex supermanifolds in the sense of Berezin and Leites; for the basic definitions and conventions, see, for example, \cite{Ber,Leites}. Let $\mathcal M=(M,\mathcal O_{\mathcal M})$ be such a supermanifold, with local even and odd coordinates $z=(z^i)$. We denote by $z^{(k)i}$, $k\geq 0$, the induced coordinates on the infinite jet bundle $J^\infty(\mathcal M)$, where
$$
z^{(0)i}=z^i,\qquad z^{(1)i}=\dot z^i.
$$
Thus, $z^{(0)i}$ are the base coordinates and $z^{(1)i}=\dot z^i$ are the corresponding fiber coordinates on the cotangent bundle $T\mathcal M$. For $k\geq2$, the coordinates $z^{(k)i}$ are the higher jet coordinates. Each $z^{(k)i}$ has the same parity as $z^i$. We refer to \cite{B} for a detailed description of $J^\infty(\mathcal M)$. If $t^i=F^i(z)$ is a change of local coordinates, then the corresponding jet coordinates are determined recursively by differentiation with respect to the even curve parameter. In particular,
\[
t^{(1)i}
=
z^{(1)a}\frac{\partial F^i}{\partial z^a}
\]
and
\[
t^{(2)i}
=
z^{(2)a}\frac{\partial F^i}{\partial z^a}
+
z^{(1)a}z^{(1)b}
\frac{\partial^2F^i}{\partial z^b\partial z^a}.
\]
Higher-order jet coordinates are obtained by repeated differentiation with respect to the even curve parameter. The important observation is that, although the coordinates $z^i$ may be even or odd, no additional signs arise from this differentiation, since the parameter on the curves used in the definition of jets is even. The products themselves are, of course, taken in the graded-commutative structure sheaf, so signs may appear if homogeneous factors are reordered.

\subsection{Formal and normalized formal exponential maps}\label{sec:formal-exponential-maps} 

The notion of a formal exponential map used here is closely related to that introduced for ordinary manifolds in \cite[Section~3.1]{CF} (see also the footnote therein for a version with a relaxed normalization condition). We adapt this notion to supermanifolds, in both the smooth and holomorphic settings.

 If $z=(z^a)$ are local even and odd coordinates on $\mathcal M$, we denote by $(z^a,\dot z^a)$ the corresponding induced coordinates on $T\mathcal M$. We consider the sheaf $\mathcal O_{\mathcal M}[[\dot z]]$ on $M$,
obtained by completing the sheaf of fiberwise polynomial functions
with respect to the natural nonnegative integer weight grading
\[
\operatorname{wt}(z^a)=0,\qquad
\operatorname{wt}(\dot z^a)=1.
\]
Thus, on an open coordinate neighborhood $U\subset M$, its sections
are formal sums
\[
F(z,\dot z)=\sum_{k=0}^{\infty}F_k(z,\dot z),
\qquad
F_k(z,\dot z)=\sum_{|I|=k}f_I(z)\dot z^I,
\]
where each $F_k$ is homogeneous of weight $k$ and
$f_I(z)\in\mathcal O_{\mathcal M}(U)$ is smooth
(respectively holomorphic) in the base coordinates.
Equivalently, the series are formal in the even fiber coordinates,
while their dependence on the odd fiber coordinates is polynomial.
In particular,
\[
\mathcal O_{\mathcal M}[[\dot z]]
=\prod_{k\geq 0}\mathcal O_{\mathcal M}[[\dot z]]_k,
\]
where $\mathcal O_{\mathcal M}[[\dot z]]_k$ denotes the
subsheaf of homogeneous elements of weight $k$.
Thus $(M,\mathcal O_{\mathcal M}[[\dot z]])$ is a formal
graded manifold, with the structure sheaf completed with respect
to the fiber weight grading.
We use the following convention for morphisms with values in
the formal sheaf $\mathcal O_{\mathcal M}[[\dot z]]$.

\begin{definition}
A morphism
$F:(M,\mathcal O_{\mathcal M}[[\dot z]])\longrightarrow(M,\mathcal O_{\mathcal M})$
is a morphism of locally ringed spaces given by a continuous map $F_0:M\to M$ and a morphism of sheaves of superalgebras
$F^\sharp:\mathcal O_{\mathcal M}\longrightarrow F_{0*}\mathcal O_{\mathcal M}[[\dot z]]$
which locally satisfies
$$
F^\sharp(f(z))=f\bigl(F^\sharp(z)\bigr)
$$
for every local superfunction $f=f(z)$, where $f\bigl(F^\sharp(z)\bigr)$ denotes the formal Taylor expansion of $f$ obtained by substituting the formal coordinate expressions $F^\sharp(z^a)$ for $z^a$. In particular, such a morphism is locally determined by the images $F^\sharp(z^a)$ of the coordinate functions.
\end{definition}
\begin{remark}
The formal substitution $f\bigl(F^\sharp(z)\bigr)$ is well defined. Indeed, setting
$F_{(0)}^\sharp:=j^\sharp\circ F^\sharp$, where 
$$j^\sharp:\mathcal O_{\mathcal M}[[\dot z]]\longrightarrow\mathcal O_{\mathcal M}$$ 
is the morphism corresponding to the zero section, locally given by setting $\dot z=0$,
we can write
$$
F^\sharp(z^a)=F_{(0)}^\sharp(z^a)+h^a,
$$
where each $h^a$ has positive fiber weight.  We define $f\bigl(F^\sharp(z)\bigr)$ by formally expanding around $F_{(0)}^\sharp(z)$ in the variables $h^a$. Since the $h^a$ have positive fiber weight, only finitely many terms contribute to each fixed weight, and hence the expansion defines an element of $\mathcal O_{\mathcal M}[[\dot z]]$. This construction applies in both the smooth and holomorphic settings; no convergence in the fiber variables is required.
\end{remark}

\begin{definition}
A \emph{formal exponential map} on $\mathcal M$ is a morphism 
$F:(M,\mathcal O_{\mathcal M}[[\dot z]])\longrightarrow(M,\mathcal O_{\mathcal M})$
whose underlying map is $\operatorname{id}_M$ and 
$F^\sharp:\mathcal O_{\mathcal M}\longrightarrow\mathcal O_{\mathcal M}[[\dot z]]$
satisfies
$j^\sharp\circ F^\sharp=\operatorname{id}_{\mathcal O_{\mathcal M}}$.
It is further required that the differential of $F$ in the fiber direction along the zero section, equivalently, the weight-one part of $F^\sharp$ in the fiber variables, be an isomorphism.
\end{definition}

In local coordinates,
$$
F^\sharp(z^a)=z^a+A^a_b(z)\dot z^b+O(\dot z^2),
$$
where $O(\dot z^2)$ denotes terms of fiber weight at least $2$ and $A=(A^a_b)$ is an invertible even matrix.

\subsection{Ehresmann and linear connections}\label{sec:Ehresmann-linear-connections} Let $\mathcal N:=T\mathcal M$ and let $p_0:\mathcal N\to\mathcal M$ denote the natural projection. We write $\mathcal T_{\mathcal M}:=\operatorname{Der}(\mathcal O_{\mathcal M})$ and $\mathcal T_{\mathcal N}:=\operatorname{Der}(\mathcal O_{\mathcal N})$ for the tangent sheaves of $\mathcal M$ and $\mathcal N$, respectively. The differential of $p_0$ gives rise to the canonical exact sequence
\[
0\longrightarrow \mathcal T_{\mathcal N}^{V}
\longrightarrow \mathcal T_{\mathcal N}
\xrightarrow{\pi}
p_0^*\mathcal T_{\mathcal M}
\longrightarrow 0,
\]
where
$p_0^*\mathcal T_{\mathcal M}
=
\mathcal O_{\mathcal N}
\otimes_{p_0^{-1}\mathcal O_{\mathcal M}}
p_0^{-1}\mathcal T_{\mathcal M}$,
$\pi$ is the morphism induced by the differential of $p_0$, and $\mathcal T_{\mathcal N}^{V}:=\ker\pi$ is the sheaf of vertical vector fields.

\begin{definition}
   An \emph{Ehresmann connection} on $T\mathcal M$ is an $\mathcal O_{\mathcal N}$-linear splitting of this sequence, that is, an $\mathcal O_{\mathcal N}$-linear morphism
$\nabla:p_0^*\mathcal T_{\mathcal M}\longrightarrow\mathcal T_{\mathcal N}$
such that $\pi\circ\nabla=\operatorname{id}$. Its image $\mathcal T_{\mathcal N}^{H}:=\nabla(p_0^*\mathcal T_{\mathcal M})$ is called the \emph{horizontal distribution} corresponding to $\nabla$. 
\end{definition}
 Thus, each Ehresmann connection gives rise to the direct sum decomposition
\[
\mathcal T_{\mathcal N}
=
\mathcal T_{\mathcal N}^{V}
\oplus
\mathcal T_{\mathcal N}^{H}.
\]
 Locally, an Ehresmann connection is determined by
$$
\nabla(\partial_a)=\partial_a-\Gamma_a^c(z,\dot z)\frac{\partial}{\partial\dot z^c},\qquad \partial_a:=\frac{\partial}{\partial z^a},
$$
for suitable functions $\Gamma_a^c\in\mathcal O_{\mathcal N}$.
An Ehresmann connection is called \emph{linear} if its coefficients are linear in the fiber coordinates, 
$$
\Gamma_a^c(z,\dot z)=\Gamma_{ab}^c(z)\dot z^b,
$$
so that the horizontal lift has the form
$$
\nabla(\partial_a)=\partial_a-\Gamma_{ab}^c(z)\dot z^b\frac{\partial}{\partial\dot z^c}.
$$

\subsection{The even operator associated with an Ehresmann connection}\label{sec Ehresmann connection} 
Each Ehresmann connection determines a decomposition $\mathcal T_{\mathcal N}=\mathcal T_{\mathcal N}^{V}\oplus\mathcal T_{\mathcal N}^{H}$, and hence the corresponding decomposition of the cotangent sheaf
$$
\Omega^1_{\mathcal N}=(\Omega^1_{\mathcal N})^H\oplus(\Omega^1_{\mathcal N})^V
$$
into its horizontal and vertical components. 
The \emph{horizontal differential}
$$
D_H:\mathcal O_{\mathcal N}\longrightarrow(\Omega^1_{\mathcal N})^H
$$
is defined as the composition of the de Rham differential $d$ with the natural projection onto the horizontal component, i.e. $D_H=\operatorname{pr}_H\circ d$,
where
$\operatorname{pr}_H:\Omega^1_{\mathcal N}\longrightarrow(\Omega^1_{\mathcal N})^H$
is the projection induced by the decomposition 
$\Omega^1_{\mathcal N}=(\Omega^1_{\mathcal N})^H\oplus(\Omega^1_{\mathcal N})^V$.

 If
$X_a:=\nabla(\partial/\partial z^a)$, then 
$D_H=dz^aX_a$.
We define the even operator associated with $\nabla$ by
$$\nabla^*:=\dot z^aX_a.
$$
This definition is independent of the choice of local coordinates. Indeed, if $t^i=F^i(z)$ is a change of coordinates, then, by the convention of Section~\ref{sec:jet-coordinates}, we have
$\dot t^i=\dot z^a\frac{\partial F^i}{\partial z^a}$.
On the other hand, we have \cite[Section 2.2.4, Corollary]{Leites}
$$
\frac{\partial}{\partial z^a}
=
\frac{\partial F^i}{\partial z^a}
\frac{\partial}{\partial t^i}.
$$
Since $\nabla$ is $\mathcal O_{\mathcal N}$-linear, we have
$$
\dot z^aX_a
=
\dot z^a\frac{\partial F^i}{\partial z^a}
\nabla\left(\frac{\partial}{\partial t^i}\right)
=
\dot t^i\nabla\left(\frac{\partial}{\partial t^i}\right).
$$
Hence, $\nabla^*$ is  a well-defined vector field on $\mathcal N$.  Since $|\dot z^a|=|X_a|=|z^a|$, every summand $\dot z^aX_a$ is even. Thus $\nabla^*$ is an even derivation of $\mathcal O_{\mathcal N}$. In particular,
$$
\nabla^*(fg)=\nabla^*(f)g+f\nabla^*(g),\qquad f,g\in\mathcal O_{\mathcal N}.
$$
If the horizontal lifts are written locally in the form
$$
X_a=\frac{\partial}{\partial z^a}-\Gamma_a^{c}(z,\dot z)\frac{\partial}{\partial\dot z^c},
$$
then
\begin{equation}\label{eq nabla^8 in coordinates}
\nabla^*=\dot z^a\frac{\partial}{\partial z^a}-\dot z^a\Gamma_a^{c}(z,\dot z)\frac{\partial}{\partial\dot z^c}.
\end{equation}
Here the order of the factors is important, since we work in the super setting. If $f=f(z)$ depends only on the base coordinates, then
\begin{equation}\label{eq nabla star on base functions}
\nabla^*(f)=\dot z^a\frac{\partial f}{\partial z^a}.
\end{equation}
According to Section~\ref{sec:jet-coordinates}, the fiber coordinates $\dot z^a$ are the first jet coordinates, $\dot z^a=z^{(1)a}$. In this notation,
$$
\nabla^*(z^{(1)c})=z^{(1)a}X_a(z^{(1)c})=-z^{(1)a}\Gamma_a^{c}(z,z^{(1)}).
$$
Thus, $\nabla^*(z^{(1)c})$ is a smooth (respectively holomorphic) function of the base coordinates and the first jet coordinates only; no higher jet coordinates occur.

 In every local coordinate chart $(z^a)$ on $\mathcal M$, with the corresponding jet coordinates $(z^{(n)a})_{n\geq0}$, define  the \emph{total jet derivative} by
\begin{equation}\label{eq total derivative}
\partial=\sum_{n=0}^{\infty}z^{(n+1)a}\frac{\partial}{\partial z^{(n)a}}.
\end{equation}
 We will show that the operators defined by \eqref{eq total derivative} in different coordinate charts agree on their overlaps and hence define a globally well-defined operator on the infinite jet space.

Let $w=w(z)$ be another system of local coordinates on $\mathcal M$, and let $w^{(k)b}$ denote the corresponding jet coordinates. Denote by $\partial_z$ and $\partial_w$ the operators defined by \eqref{eq total derivative} in the $z$- and $w$-coordinates, respectively. Thus,
$$
\partial_z
=
\sum_{n=0}^{\infty}
z^{(n+1)a}\frac{\partial}{\partial z^{(n)a}},
\qquad
\partial_w
=
\sum_{k=0}^{\infty}
w^{(k+1)b}\frac{\partial}{\partial w^{(k)b}}.
$$
We prove that $\partial_z=\partial_w$ on the overlap.

By the definition of jet coordinates in terms of curves with an even parameter $t$, for any such curve we have
$$
z^{(n+1)a}=\frac{d}{dt}z^{(n)a},
\qquad
w^{(k+1)b}=\frac{d}{dt}w^{(k)b}.
$$
Here and below, these identities are understood along an arbitrary curve defining the corresponding jet.
Since $w^{(k)b}$ depends only on $z^{(0)},\ldots,z^{(k)}$, differentiation along the curve gives
$$
\frac{d}{dt}w^{(k)b}
=
\sum_{n=0}^{k}
\frac{d z^{(n)a}}{dt}
\frac{\partial w^{(k)b}}{\partial z^{(n)a}}
=
\sum_{n=0}^{k}
z^{(n+1)a}
\frac{\partial w^{(k)b}}{\partial z^{(n)a}}.
$$
Since the curve is arbitrary, this gives the following relation between the jet coordinates:
$$
w^{(k+1)b}
=
\sum_{n=0}^{k}
z^{(n+1)a}
\frac{\partial w^{(k)b}}{\partial z^{(n)a}}.
$$
On the other hand, by the definition of $\partial_z$,
$$
\partial_z\bigl(w^{(k)b}\bigr)
=
\sum_{n=0}^{k}
z^{(n+1)a}
\frac{\partial w^{(k)b}}{\partial z^{(n)a}}
=
w^{(k+1)b} = \partial_w\bigl(w^{(k)b}\bigr).
$$
Thus the operators defined by \eqref{eq total derivative} in different coordinate charts agree on overlaps, and hence $\partial$ is a globally well-defined even derivation on the infinite jet space.

\section{Linearization of jets} 

In this section, we prove a linearization theorem for the infinite jet
space for supermanifolds. We show that a choice of an Ehresmann
connection allows us to construct new jet coordinates which transform
linearly under changes of local coordinates. As a result, we obtain an
isomorphism between the infinite jet space and its linear model
\[
J^\infty_{\mathrm{lin}}(\mathcal M)
=
\bigoplus_{n=0}^{\infty}T[n]\mathcal M.
\]

Let us start with an example.

\subsection{Linearization of second jets}
\label{ex k=2 linearization for jets}
Let $t^i=F^i(z)$ be a change of local even and odd coordinates on $\mathcal M$. According to the convention for jet coordinates fixed in Section~\ref{sec:jet-coordinates}, the corresponding first jet coordinates transform as
$$
t^{(1)i}
=
z^{(1)a}\frac{\partial F^i}{\partial z^a}.
$$
Applying the total jet derivative  \eqref{eq total derivative}  to this transformation law, we obtain
$$
t^{(2)i}
=
z^{(2)a}\frac{\partial F^i}{\partial z^a}
+
z^{(1)a}z^{(1)b}
\frac{\partial^2F^i}{\partial z^b\partial z^a}.
$$
We define the modified second jet coordinates by
$$
\widetilde z^{(2)a}
=
z^{(2)a}-\nabla^*(z^{(1)a}),
\qquad
\widetilde t^{(2)i}
=
t^{(2)i}-\nabla^*(t^{(1)i}),
$$
where $\nabla^*$ is the even operator associated with the Ehresmann connection; see Section~\ref{sec Ehresmann connection}. Since $\nabla^*$ is an even derivation, using the transformation law for the first jet coordinates we obtain
$$
\begin{aligned}
\nabla^*(t^{(1)i})
&=
\nabla^*\left(
z^{(1)a}\frac{\partial F^i}{\partial z^a}
\right)
\\
&=
\nabla^*(z^{(1)a})\frac{\partial F^i}{\partial z^a}
+
z^{(1)a}\nabla^*\left(
\frac{\partial F^i}{\partial z^a}
\right).
\end{aligned}
$$
Since ${\partial F^i}/{\partial z^a}$ depends only on the base coordinates, formula~(\ref{eq nabla star on base functions}) gives
$$
\nabla^*\left(
\frac{\partial F^i}{\partial z^a}
\right)
=
z^{(1)b}
\frac{\partial^2F^i}{\partial z^b\partial z^a}.
$$
Therefore,
$$
\nabla^*(t^{(1)i})
=
\nabla^*(z^{(1)a})\frac{\partial F^i}{\partial z^a}
+
z^{(1)a}z^{(1)b}
\frac{\partial^2F^i}{\partial z^b\partial z^a}.
$$
We now compute the transition law for the modified second jet coordinates:
$$
\begin{aligned}
\widetilde t^{(2)i}
&=
t^{(2)i}-\nabla^*(t^{(1)i})
\\
&=
z^{(2)a}\frac{\partial F^i}{\partial z^a}
+
z^{(1)a}z^{(1)b}
\frac{\partial^2F^i}{\partial z^b\partial z^a}
\\
&\quad
-
\nabla^*(z^{(1)a})\frac{\partial F^i}{\partial z^a}
-
z^{(1)a}z^{(1)b}
\frac{\partial^2F^i}{\partial z^b\partial z^a}
\\
&=
\left(
z^{(2)a}-\nabla^*(z^{(1)a})
\right)
\frac{\partial F^i}{\partial z^a}
=
\widetilde z^{(2)a}\frac{\partial F^i}{\partial z^a}.
\end{aligned}
$$
Thus,
$$
\widetilde t^{(2)i}
=
\widetilde z^{(2)a}\frac{\partial F^i}{\partial z^a},
$$
and the modified second jet coordinates transform linearly, with the same transition law as the fiber coordinates on $T\mathcal M$.

\subsection{Linearization of higher jets}\label{sec:linearization-higher-jets}
We now extend the construction of the modified second jet coordinates to higher orders. Suppose that, for some $k\geq2$, modified jet coordinates
$$
\widetilde z^{(1)a}=z^{(1)a},\qquad
\widetilde z^{(2)a},\ldots,\widetilde z^{(k)a}
$$
have already been constructed and transform linearly. More precisely, for every $r=1,\ldots,k$, under a change of coordinates $t^i=F^i(z)$ we have
$$
\widetilde t^{(r)i}
=
\widetilde z^{(r)a}\frac{\partial F^i}{\partial z^a}.
$$
We will construct the modified jet coordinates $\widetilde z^{(k+1)a}$ and prove that they satisfy the same transformation law.

In particular, the transformation law for $\widetilde z^{(k)a}$ is the same as that of the fiber coordinates $\dot z^a$ on $T\mathcal M$. Hence the pair
$\bigl(z^a,\widetilde z^{(k)a}\bigr)$
may be regarded as a system of coordinates on a copy of $T\mathcal M$. We can therefore apply a given Ehresmann connection to this copy. If
$$
X_a
=
\frac{\partial}{\partial z^a}
-
\Gamma_a^c(z,\dot z)
\frac{\partial}{\partial\dot z^c}
$$
are its horizontal lifts in the coordinates $(z^a,\dot z^a)$, then in the coordinates $(z^a,\widetilde z^{(k)a})$ they take the form
$$
X_a^{(k)}
=
\frac{\partial}{\partial z^a}
-
\Gamma_a^c\bigl(z,\widetilde z^{(k)}\bigr)
\frac{\partial}{\partial\widetilde z^{(k)c}}.
$$
As in the construction for the second jet coordinates, we now use the first jet coordinates $z^{(1)a}$ and define
$$
\nabla^{*(k)}
=
z^{(1)a}X_a^{(k)}.
$$
This is an even operator. Its action on the fiber coordinates $\widetilde z^{(k)c}$ is given by
$$
\nabla^{*(k)}
\bigl(\widetilde z^{(k)c}\bigr)
=
-z^{(1)a}
\Gamma_a^c\bigl(z,\widetilde z^{(k)}\bigr).
$$

\begin{theorem}\label{theor Ehresmann jet space}
Let \(\mathcal M\) be a (complex) supermanifold endowed with an Ehresmann connection as above. Then, for every \(N\), there exists a triangular change of coordinates on \(J^N(\mathcal M)\) such that the transition functions in the new coordinates are linear. These changes of coordinates are compatible with the natural projections between the finite jet spaces and therefore induce an isomorphism
$$
J^\infty(\mathcal M)\cong\bigoplus_{n=0}^{\infty}T[n]\mathcal M,
$$
where the structure sheaves on both sides are understood as the direct limits of the corresponding finite-order structure sheaves. In general, the change of coordinates is not homogeneous with respect to the jet grading.
\end{theorem}

\begin{proof}
We construct local coordinates on \(J^\infty(\mathcal M)\) in which
the transition functions become linear. Let
\((z^i,z^{(1)i},z^{(2)i},\ldots)\) be the standard jet coordinates
associated with local coordinates \((z^i)\) on \(\mathcal M\).

\smallskip

{\bf Base case.}
The base case was established in
Section~\ref{ex k=2 linearization for jets}. More precisely, setting
$\widetilde z^{(1)a}=z^{(1)a}$,
we define
\[
\widetilde z^{(2)a}
=
\partial\bigl(\widetilde z^{(1)a}\bigr)
-
\nabla^*\bigl(\widetilde z^{(1)a}\bigr)
=
z^{(2)a}-\nabla^*(z^{(1)a}).
\]
As shown in Section~\ref{ex k=2 linearization for jets}, these modified
second jet coordinates transform linearly:
\[
\widetilde t^{(2)i}
=
\frac{\partial F^i}{\partial z^a}\,
\widetilde z^{(2)a}.
\]
Hence $\widetilde z^{(2)a}$ transform as fiber coordinates of $T\mathcal M$. Applying the construction described above for $k=2$, we obtain the corresponding horizontal lifts $X_a^{(2)}$ and the even operator
$\nabla^{*(2)}=z^{(1)a}X_a^{(2)}$.

\smallskip

{\bf Induction hypothesis.}
Assume that the modified jet coordinates
$\widetilde z^{(k)a}$, $1\leq k\leq n-1$, have been constructed and transform linearly:
$$
\widetilde t^{(k)i}
=
\widetilde z^{(k)a}\frac{\partial F^i}{\partial z^a},
\qquad
1\leq k\leq n-1.
$$

\smallskip
{\bf Inductive step.}
Now $\widetilde z^{(n-1)a}$ transform in the same way as the fiber coordinates on $T\mathcal M$. Hence the pair
$\bigl(z^a,\widetilde z^{(n-1)a}\bigr)$
may be regarded as coordinates on a copy of $T\mathcal M$. Using the given Ehresmann connection $\nabla$ on $T\mathcal M$, we obtain on this copy the horizontal lifts
$$
X_a^{(n-1)}
=
\frac{\partial}{\partial z^a}
-
\Gamma_a^c\bigl(z,\widetilde z^{(n-1)}\bigr)
\frac{\partial}{\partial\widetilde z^{(n-1)c}},
$$
and define
$\nabla^{*(n-1)}
=
z^{(1)a}X_a^{(n-1)}$.
Consider the transformation law for the $(n-1)$-st modified jet coordinates,
$$
\widetilde t^{(n-1)i}
=
\widetilde z^{(n-1)a}\frac{\partial F^i}{\partial z^a}.
$$
Applying the total derivative $\partial$ to this equality, we obtain
$$
\partial\bigl(\widetilde t^{(n-1)i}\bigr)
=
\partial\bigl(\widetilde z^{(n-1)a}\bigr)
\frac{\partial F^i}{\partial z^a}
+
\widetilde z^{(n-1)a}z^{(1)b}
\frac{\partial^2F^i}{\partial z^b\partial z^a}.
$$
On the other hand, applying $\nabla^{*(n-1)}$ to the same transformation law and using the fact that $\nabla^{*(n-1)}$ is even, we obtain
$$
\begin{aligned}
\nabla^{*(n-1)}
\bigl(\widetilde t^{(n-1)i}\bigr)
&=
\nabla^{*(n-1)}
\left(
\widetilde z^{(n-1)a}
\frac{\partial F^i}{\partial z^a}
\right)
\\
&=
\nabla^{*(n-1)}
\bigl(\widetilde z^{(n-1)a}\bigr)
\frac{\partial F^i}{\partial z^a}
+
\widetilde z^{(n-1)a}
\nabla^{*(n-1)}
\left(
\frac{\partial F^i}{\partial z^a}
\right).
\end{aligned}
$$
Since ${\partial F^i}/{\partial z^a}$ depends only on the base coordinates,
$$
\nabla^{*(n-1)}
\left(
\frac{\partial F^i}{\partial z^a}
\right)
=
z^{(1)b}
\frac{\partial^2F^i}{\partial z^b\partial z^a}.
$$
Therefore,
$$
\begin{aligned}
\nabla^{*(n-1)}
\bigl(\widetilde t^{(n-1)i}\bigr)
&=
\nabla^{*(n-1)}
\bigl(\widetilde z^{(n-1)a}\bigr)
\frac{\partial F^i}{\partial z^a}
+
\widetilde z^{(n-1)a}z^{(1)b}
\frac{\partial^2F^i}{\partial z^b\partial z^a}.
\end{aligned}
$$
We now define the modified $n$-th jet coordinates by
$$
\widetilde z^{(n)a}
=
\partial\bigl(\widetilde z^{(n-1)a}\bigr)
-
\nabla^{*(n-1)}
\bigl(\widetilde z^{(n-1)a}\bigr),
$$
and similarly
$$
\widetilde t^{(n)i}
=
\partial\bigl(\widetilde t^{(n-1)i}\bigr)
-
\nabla^{*(n-1)}
\bigl(\widetilde t^{(n-1)i}\bigr).
$$
Substituting the two formulas above, we obtain
$$
\begin{aligned}
\widetilde t^{(n)i}
&=
\partial\bigl(\widetilde z^{(n-1)a}\bigr)
\frac{\partial F^i}{\partial z^a}
+
\widetilde z^{(n-1)a}z^{(1)b}
\frac{\partial^2F^i}{\partial z^b\partial z^a}
\\
&\quad
-
\nabla^{*(n-1)}
\bigl(\widetilde z^{(n-1)a}\bigr)
\frac{\partial F^i}{\partial z^a}
-
\widetilde z^{(n-1)a}z^{(1)b}
\frac{\partial^2F^i}{\partial z^b\partial z^a}
\\
&=
\left(
\partial\bigl(\widetilde z^{(n-1)a}\bigr)
-
\nabla^{*(n-1)}
\bigl(\widetilde z^{(n-1)a}\bigr)
\right)
\frac{\partial F^i}{\partial z^a}
=
\widetilde z^{(n)a}\frac{\partial F^i}{\partial z^a}.
\end{aligned}
$$
The terms containing the second derivatives of $F^i$ occur in exactly the same order and therefore cancel directly. No permutation of homogeneous factors is involved, and hence no additional Koszul sign appears. Thus the modified $n$-th jet coordinates transform linearly:
$$
\widetilde t^{(n)i}
=
\widetilde z^{(n)a}\frac{\partial F^i}{\partial z^a}.
$$
Consequently, $\widetilde z^{(n)a}$ again transform as fiber coordinates on $T\mathcal M$, so the same construction can be applied at the next step. This completes the induction.

\smallskip

{\bf Triangularity.} It remains to verify that the resulting change of jet coordinates is triangular. By the induction hypothesis,
$$
\widetilde z^{(n-1)a}
=
z^{(n-1)a}+P^{(n-1)a},
$$
where $P^{(n-1)a}$ depends only on jet coordinates of order less than $n-1$. Therefore,
$$
\partial\bigl(\widetilde z^{(n-1)a}\bigr)
=
z^{(n)a}
+
\text{terms involving only jet coordinates of order at most }n-1.
$$
Moreover,
$$
\nabla^{*(n-1)}
\bigl(\widetilde z^{(n-1)a}\bigr)
=
-z^{(1)b}
\Gamma_b^a\bigl(z,\widetilde z^{(n-1)}\bigr)
$$
is a smooth, respectively holomorphic, function of jet coordinates of order at most $n-1$. Hence
$$
\widetilde z^{(n)a}
=
z^{(n)a}+P^{(n)a},
$$
where $P^{(n)a}$ depends only on jet coordinates of order at most $n-1$. Thus, the change of coordinates is triangular with identity on the highest-order term. It is therefore invertible  and defines a smooth (respectively holomorphic) change of coordinates on each finite jet space.
\end{proof}


\begin{remark}
The construction above is performed on each finite jet space $J^N(\mathcal M)$. Since the triangular changes of coordinates are compatible with the natural projections
$J^{N+1}(\mathcal M)\longrightarrow J^N(\mathcal M)$,
they induce the isomorphism of infinite jet spaces stated in Theorem~\ref{theor Ehresmann jet space}.

More precisely, the infinite jet space $J^\infty(\mathcal M)$ is the inverse limit of the finite jet spaces $J^N(\mathcal M)$. Its structure sheaf is obtained as the direct limit of the corresponding structure sheaves, pulled back to the inverse-limit space:
$$
\mathcal O_{J^\infty(\mathcal M)}
=
\varinjlim_N\mathcal O_{J^N(\mathcal M)}.
$$
Here, the pullbacks along the natural projections are understood. In particular, every local section depends on jet coordinates of finite order.
\end{remark}

\begin{remark}\label{rem:unrestricted-formal-completion}
The triangular change of coordinates does not, in general, extend to the algebra of formal power series in all jet coordinates. To illustrate this, consider the one-dimensional manifold with coordinate $x$ and the Ehresmann connection
$$
X=\frac{\partial}{\partial x}
-\Gamma(x,\dot x)\frac{\partial}{\partial\dot x},
\qquad \Gamma(x,\dot x)=1.
$$
By the construction in the proof of Theorem~\ref{theor Ehresmann jet space}, the modified jet coordinates satisfy
$$
\widetilde x^{(1)}=x^{(1)},\qquad
\widetilde x^{(n)}
=
\partial\bigl(\widetilde x^{(n-1)}\bigr)
-
\nabla^{*(n-1)}
\bigl(\widetilde x^{(n-1)}\bigr),
\qquad n\geq 2.
$$
Since $\Gamma(x,\dot x)=1$, we have
$\nabla^{*(n-1)}
\bigl(\widetilde x^{(n-1)}\bigr)
=
-x^{(1)}$.
Consequently,
$$
\widetilde x^{(n)}
=
\partial\bigl(\widetilde x^{(n-1)}\bigr)+x^{(1)}.
$$
Using the definition of the total jet derivative \eqref{eq total derivative},  we obtain
$$
\widetilde x^{(1)}=x^{(1)},\qquad
\widetilde x^{(2)}=x^{(2)}+x^{(1)},\qquad
\widetilde x^{(3)}=x^{(3)}+x^{(2)}+x^{(1)},
$$
and, by induction,
$$
\widetilde x^{(n)}
=
x^{(n)}+x^{(n-1)}+\cdots+x^{(1)}.
$$
At every finite jet order, this defines an invertible triangular change of coordinates. However, consider the formal series
$$
F=\sum_{n=1}^{\infty}\widetilde x^{(n)}.
$$
Substituting the expressions above, we obtain
$$
F
=
\sum_{n=1}^{\infty}\sum_{k=1}^{n}x^{(k)}.
$$
In particular, the coefficient of $x^{(1)}$ would be the infinite sum
$$
1+1+1+\cdots,
$$
which is not defined in the formal power series algebra. Consequently, although the triangular transformations are well-defined and invertible at every finite jet order, they do not necessarily extend to an isomorphism of the unrestricted formal completions in all jet variables.
\end{remark}

\subsection{Formal completion of infinite jet spaces}
\label{sec:formal-completion-jets}

We now introduce the sheaves of formal power series in jet
coordinates. Assign weight $k$ to the jet coordinates
$z^{(k)a}$, $k\geq 1$, and weight zero to functions
on $\mathcal M$. Let $\mathcal J^\infty_{\mathcal M}$
denote the completion of the sheaf of polynomial jet
superalgebras with respect to the total jet weight.
Thus, its local sections are formal sums
\[
f=\sum_{m=0}^{\infty}f_m,
\]
where $f_m$ is a homogeneous polynomial of total jet
weight $m$, with coefficients in $\mathcal O_{\mathcal M}$.
Similarly, let $\mathcal J^\infty_{\mathrm{lin},\mathcal M}$
denote the completion of the sheaf of polynomial
superalgebras in the linearized jet coordinates
$\widetilde z^{(k)a}$, assigning weight $k$ to
$\widetilde z^{(k)a}$.

\begin{theorem}\label{theor Ehresmann completed jet space}
Suppose that the Ehresmann connection preserves the zero section, i.e.,
$\Gamma_a^c(z,0)=0$.
Then the triangular change of coordinates constructed in Theorem~\ref{theor Ehresmann jet space} extends uniquely to an isomorphism of sheaves of formal superalgebras on $\mathcal M$,
$$
\mathcal J^\infty_{\mathcal M}
\simeq
\mathcal J^\infty_{\mathrm{lin},\mathcal M}.
$$
\end{theorem}

\begin{proof}
Recall that the triangular change of coordinates constructed in the proof of Theorem~\ref{theor Ehresmann jet space} is given recursively by
$$
\widetilde z^{(1)a}=z^{(1)a},\qquad
\widetilde z^{(k+1)a}
=
\partial\bigl(\widetilde z^{(k)a}\bigr)
-
\nabla^{*(k)}\bigl(\widetilde z^{(k)a}\bigr).
$$
We show by induction that
$$
\widetilde z^{(k)a}
=
z^{(k)a}
+
\text{terms of jet weight at least }k.
$$
The assertion is clear for $k=1$. Suppose that it holds for $k$. Since the total jet derivative $\partial$, defined in \eqref{eq total derivative}, raises the jet weight by one, we obtain
$$
\partial\bigl(\widetilde z^{(k)a}\bigr)
=
z^{(k+1)a}
+
\text{terms of jet weight at least }k+1.
$$
On the other hand, by the definition of $\nabla^{*(k)}$ in the proof of Theorem~\ref{theor Ehresmann jet space},
$$
\nabla^{*(k)}\bigl(\widetilde z^{(k)c}\bigr)
=
-z^{(1)a}\Gamma_a^c\bigl(z,\widetilde z^{(k)}\bigr).
$$
Since the connection preserves the zero section, $\Gamma_a^c(z,0)=0$, the Taylor expansion of $\Gamma_a^c(z,\widetilde z^{(k)})$ contains only terms of positive degree in $\widetilde z^{(k)}$. By the induction hypothesis, these terms have jet weight at least $k$. Consequently, $\nabla^{*(k)}\bigl(\widetilde z^{(k)c}\bigr)$ has jet weight at least $k+1$. It follows that
$$
\widetilde z^{(k+1)a}
=
z^{(k+1)a}
+
\text{terms of jet weight at least }k+1.
$$
Thus, the triangular change of coordinates preserves the jet-weight filtration.

The inverse triangular change of coordinates exists by Theorem~\ref{theor Ehresmann jet space}, and it satisfies the same weight condition:
$$
z^{(k)a}
=
\widetilde z^{(k)a}
+
\text{terms of jet weight at least }k.
$$
Consequently, when either the coordinate transformation or its inverse is substituted into a formal series, only finitely many terms can contribute to any fixed jet weight. Thus both substitutions are well-defined on the sheaves of formal superalgebras. Since they are mutually inverse at every finite jet order, their extensions are also mutually inverse. This proves the theorem.
\end{proof}

\begin{corollary}\label{cor theorem linear}
Suppose that the Ehresmann connection in Theorem~\ref{theor Ehresmann jet space} is induced by a linear connection on $\mathcal M$. Then the triangular change of jet coordinates is homogeneous with respect to the jet grading and extends to an isomorphism of sheaves of formal superalgebras
$$
\mathcal J^\infty_{\mathcal M}
\simeq
\mathcal J^\infty_{\mathrm{lin},\mathcal M}
$$
preserving the jet grading.
\end{corollary}

\begin{proof}
Since the connection is linear, its coefficients have the form
$$
\Gamma_a^c(z,\dot z)=\Gamma_{ab}^c(z)\dot z^b.
$$
In particular, $\Gamma_a^c(z,0)=0$, so the extension to the sheaves of formal superalgebras follows from Theorem~\ref{theor Ehresmann completed jet space}.

It remains to prove homogeneity. We proceed by induction on $k$. The first modified jet coordinates satisfy $\widetilde z^{(1)a}=z^{(1)a}$ and are therefore homogeneous of jet degree $1$. Suppose that $\widetilde z^{(k)a}$ are homogeneous of jet degree $k$. Then
$$
\nabla^{*(k)}
\bigl(\widetilde z^{(k)c}\bigr)
=
-z^{(1)a}\Gamma_{ab}^c(z)\widetilde z^{(k)b}
$$
is homogeneous of jet degree $k+1$. Since the total jet derivative $\partial$, defined in \eqref{eq total derivative}, increases the jet degree by one, the recursive formula
$$
\widetilde z^{(k+1)c}
=
\partial\bigl(\widetilde z^{(k)c}\bigr)
-
\nabla^{*(k)}
\bigl(\widetilde z^{(k)c}\bigr)
$$
shows that $\widetilde z^{(k+1)c}$ is homogeneous of jet degree $k+1$. This completes the induction.
\end{proof}

\section{Covariant total derivative}

\subsection{The covariant total derivative}\label{sec:covariant-total-derivative}
Recall that, in the proof of Theorem~\ref{theor Ehresmann jet space}, each set of modified coordinates $\widetilde z^{(k)a}$ was identified with the fiber coordinates of a copy of $T\mathcal M$. On each such copy, the given Ehresmann connection determines the horizontal lifts
$$
X_a^{(k)}
=
\frac{\partial}{\partial z^a}
-
\Gamma_a^c\bigl(z,\widetilde z^{(k)}\bigr)
\frac{\partial}{\partial\widetilde z^{(k)c}}.
$$
Together, these horizontal lifts determine an even derivation on the linearized infinite jet space, which we continue to denote by $\nabla^*$. Its local expression is
$$
\nabla^*
=
z^{(1)a}\frac{\partial}{\partial z^a}
-
\sum_{k\geq1}
z^{(1)a}
\Gamma_a^c\bigl(z,\widetilde z^{(k)}\bigr)
\frac{\partial}{\partial\widetilde z^{(k)c}}.
$$
The infinite sum is well-defined on the structure sheaf of the infinite jet space, since every local section depends on jet coordinates of finite order and therefore only finitely many terms contribute. Equivalently, $\nabla^*$ is determined by
$$
\nabla^*(z^c)=z^{(1)c},
\qquad
\nabla^*(\widetilde z^{(k)c})
=
-z^{(1)a}\Gamma_a^c\bigl(z,\widetilde z^{(k)}\bigr),
\qquad k\geq1,
$$
and the Leibniz rule. In particular, for a function $f=f(z)$ depending only on the base coordinates, we have
$$
\nabla^*(f)=z^{(1)a}\frac{\partial f}{\partial z^a},
$$
in agreement with \eqref{eq nabla star on base functions}.

Using the triangular change of coordinates of Theorem~\ref{theor Ehresmann jet space}, we regard $\nabla^*$ as an operator on the original infinite jet space. We define the \emph{covariant total derivative} associated with the Ehresmann connection by
\begin{equation}\label{eq covariant total derivative}
\partial_\nabla:=\partial-\nabla^*,
\end{equation}
where $\partial$ is the total jet derivative defined in \eqref{eq total derivative}. Since both operators are even, $\partial_\nabla$ is an even derivation. The recursive construction in the proof of Theorem~\ref{theor Ehresmann jet space} can now be written as
$$
\widetilde z^{(1)a}=z^{(1)a},
\qquad
\widetilde z^{(k+1)a}
=
\partial_\nabla\bigl(\widetilde z^{(k)a}\bigr),
\qquad k\geq1.
$$
Consequently,
$$
\widetilde z^{(n)a}
=
\partial_\nabla^{\,n-1}\bigl(z^{(1)a}\bigr),
\qquad n\geq1.
$$
Thus the isomorphism constructed in Theorem~\ref{theor Ehresmann jet space} is explicitly described by
\begin{equation}\label{eq jet linearization map}
\begin{gathered}
\Phi:
J^\infty(\mathcal M)
\longrightarrow
J^\infty_{\mathrm{lin}}(\mathcal M)
=
\bigoplus_{n=0}^{\infty}T[n]\mathcal M,
\\[2mm]
\widetilde z^{(0)a}=z^a,
\qquad
\widetilde z^{(1)a}=z^{(1)a},
\qquad
\widetilde z^{(n)a}
=
\partial_\nabla^{\,n-1}\bigl(z^{(1)a}\bigr),
\quad n\geq2.
\end{gathered}
\end{equation}

For a general Ehresmann connection, the coefficients $\Gamma_a^c(z,\widetilde z^{(k)})$ may depend nonlinearly on the fiber coordinates. Consequently, $\nabla^*$ need not be homogeneous with respect to the jet grading, and $\partial_\nabla$ does not necessarily increase the jet weight by exactly one. The resulting change of coordinates is therefore triangular but need not be graded homogeneous. If the Ehresmann connection is induced by a linear connection, the change of coordinates is homogeneous with respect to the jet grading, as established in Corollary~\ref{cor theorem linear}.

\subsection{The ideal $\mathcal A$}\label{sec:ideal-A}
Let $\mathcal A\subset\mathcal O_{J^\infty_{\mathrm{lin}}(\mathcal M)}$ be the ideal sheaf locally generated by the linearized jet coordinates $\widetilde z^{(k)a}$, $k>1$. Thus
$$
\mathcal A=\left\langle \widetilde z^{(k)a}\,;\ k>1\right\rangle.
$$
Since the modified jet coordinates transform linearly, this ideal is globally defined. The quotient by $\mathcal A$ is naturally identified with the structure sheaf of $T\mathcal M$:
$$
\mathcal O_{J^\infty_{\mathrm{lin}}(\mathcal M)}/\mathcal A
\simeq
\mathcal O_{T\mathcal M}.
$$

\begin{theorem}\label{theor ideal A}
Let $\mathcal M$ be endowed with an Ehresmann connection, and let $\mathcal A$ be the ideal defined above.
\begin{enumerate}
\item The covariant total derivative $\partial_\nabla$ preserves $\mathcal A$ and induces the zero derivation on $\mathcal O_{T\mathcal M}$:
$$
\bar\partial_\nabla=0.
$$
\item If, in addition, the Ehresmann connection preserves the zero section, i.e., $\Gamma_a^c(z,0)=0$, then $\nabla^*$ and $\partial$ also preserve $\mathcal A$. The induced operators $\bar\nabla^*$ and $\bar\partial$  satisfy
$$
\bar\partial=\bar\nabla^*=\nabla^*,
$$
where on the right-hand side $\nabla^*$ denotes the original connection operator on $T\mathcal M$.
\end{enumerate}
\end{theorem}

\begin{proof}
By the recursive construction of the modified jet coordinates,
$\partial_\nabla\bigl(\widetilde z^{(k)a}\bigr)
=
\widetilde z^{(k+1)a}$.
Hence $\partial_\nabla(\mathcal A)\subset\mathcal A$.
Furthermore, since
$$
\partial_\nabla(z^a)=0,\qquad
\partial_\nabla(\dot z^a)=\widetilde z^{(2)a}\in\mathcal A,
$$
the induced derivation $\bar\partial_\nabla$ vanishes on the coordinate functions of the quotient, and hence $\bar\partial_\nabla=0$.
Now suppose that $\Gamma_a^c(z,0)=0$. Since
$$
\nabla^*\bigl(\widetilde z^{(k)c}\bigr)
=
-z^{(1)a}\Gamma_a^c\bigl(z,\widetilde z^{(k)}\bigr),
$$
for $k>1$, the right-hand side belongs to $\mathcal A$. Hence, $\nabla^*(\mathcal A)\subset\mathcal A$. Since $\partial=\partial_\nabla+\nabla^*$, we also have $\partial(\mathcal A)\subset\mathcal A$.
The induced operators satisfy
$$
\bar\partial
=
\bar\partial_\nabla+\bar\nabla^*
=
\bar\nabla^*.
$$

Finally, after passing to the quotient by $\mathcal A$, only the $k=1$ term in the local expression for $\nabla^*$ remains, giving
$$
\bar\nabla^*
=
\dot z^a\frac{\partial}{\partial z^a}
-
\dot z^a\Gamma_a^c(z,\dot z)
\frac{\partial}{\partial\dot z^c}.
$$
This is precisely the original connection operator $\nabla^*$ on $T\mathcal M$.
\end{proof}

\begin{corollary}\label{cor ideal A linear}
If the Ehresmann connection is induced by a linear connection on $\mathcal M$, then all the conclusions of Theorem~\ref{theor ideal A} hold.
\end{corollary}

\begin{proof}
For a linear connection, $\Gamma_a^c(z,\dot z)=\Gamma_{ab}^c(z)\dot z^b$, and therefore $\Gamma_a^c(z,0)=0$.
\end{proof}

\begin{remark}
The ideal $\mathcal A$ has a natural analogue in the sheaf of formal superalgebras $\mathcal J^\infty_{\mathrm{lin},\mathcal M}$ introduced in Section~\ref{sec:formal-completion-jets}. Let $\widehat{\mathcal A}$ be the ideal consisting of formal series whose monomials contain at least one coordinate $\widetilde z^{(k)a}$ with $k>1$. Then, locally,
$$
\mathcal J^\infty_{\mathrm{lin},\mathcal M}/\widehat{\mathcal A}
\simeq
\mathcal O_{\mathcal M}[[\dot z]],
$$
where $\dot z^a=\widetilde z^{(1)a}$, and where $\dot z^a=\widetilde z^{(1)a}$, and
$\mathcal O_{\mathcal M}[[\dot z]]$ is the sheaf of formal power series
introduced in Section~\ref{sec:formal-completion-jets}.

The conclusions of Theorem~\ref{theor ideal A} remain valid in this setting. Indeed, $\partial_\nabla$ preserves $\widehat{\mathcal A}$ and induces the zero derivation on the quotient. If, in addition, $\Gamma_a^c(z,0)=0$, then $\nabla^*$ and $\partial$ also preserve $\widehat{\mathcal A}$, and the induced operators satisfy
$$
\bar\partial_\nabla=0,
\qquad
\bar\partial=\bar\nabla^*=\nabla^*,
$$
where on the right-hand side $\nabla^*$ denotes the formal extension of the original connection operator on $T\mathcal M$ to $\mathcal O_{\mathcal M}[[\dot z]]$.
\end{remark}


\subsection{Higher jets and the connection operator}
In this subsection, we assume that the Ehresmann connection preserves the zero section, that is, $\Gamma_a^c(z,0)=0$. By Theorem~\ref{theor ideal A}, the operators $\nabla^*$, $\partial_\nabla$, and $\partial$ preserve the ideal $\mathcal A$, and the operator induced by $\partial$ on the quotient coincides with the original connection operator $\nabla^*$ on $T\mathcal M$.
Let
$$
\pi_{\mathcal A}:
\mathcal O_{J^\infty_{\mathrm{lin}}(\mathcal M)}
\longrightarrow
\mathcal O_{J^\infty_{\mathrm{lin}}(\mathcal M)}/\mathcal A
$$
be the canonical quotient map. We use the identification of the quotient with the structure sheaf of $T\mathcal M$ established in Section~\ref{sec:ideal-A}.

\begin{theorem}\label{thm inverse jets mod A}
If the Ehresmann connection preserves the zero section, then
the isomorphism  \eqref{eq jet linearization map}
satisfies
$$
\left(\pi_{\mathcal A}\circ(\Phi^\sharp)^{-1}\right)
\bigl(z^{(n)a}\bigr)
=
(\nabla^*)^n(z^a),
\qquad n\geq0,
$$
where $\nabla^*$ is the original connection operator on $T\mathcal M$.
\end{theorem}

\begin{proof}
For $n=0$, the assertion is immediate. For $n=1$, we have
$$
\left(\pi_{\mathcal A}\circ(\Phi^\sharp)^{-1}\right)
\bigl(z^{(1)a}\bigr)
=
\dot z^a
=
\nabla^*(z^a).
$$
Suppose that the assertion holds for some $n\geq1$. Since $z^{(n+1)a}=\partial(z^{(n)a})$ and $\partial$ preserves $\mathcal A$, we may pass this identity to the quotient. By Theorem~\ref{theor ideal A}, the induced operator satisfies $\bar\partial=\nabla^*$. Hence
$$
\begin{aligned}
\left(\pi_{\mathcal A}\circ(\Phi^\sharp)^{-1}\right)
\bigl(z^{(n+1)a}\bigr)
&=
\bar\partial\left(
\left(\pi_{\mathcal A}\circ(\Phi^\sharp)^{-1}\right)
\bigl(z^{(n)a}\bigr)\right)\\
&=
\nabla^*\bigl((\nabla^*)^n(z^a)\bigr)
=
(\nabla^*)^{n+1}(z^a).
\end{aligned}
$$
This completes the induction.
\end{proof}

\subsection{The formal exponential map associated with higher jets}
Recall that $\mathcal J^\infty_{\mathcal M}$  is the sheaf of formal jet series introduced in Section~\ref{sec:formal-completion-jets}. The canonical Taylor morphism is defined by
$$
\tau^\sharp:
\mathcal O_{\mathcal M}\longrightarrow
\mathcal J^\infty_{\mathcal M},
\qquad
\tau^\sharp(f)
=
\sum_{n=0}^{\infty}\frac{1}{n!}\partial^n(f).
$$
In particular, for a local coordinate $z^a$,
$$
\tau^\sharp(z^a)
=
\sum_{n=0}^{\infty}\frac{1}{n!}z^{(n)a},
\qquad z^{(0)a}=z^a.
$$
By Theorem~\ref{theor Ehresmann completed jet space}, the assumption
$\Gamma_a^c(z,0)=0$ ensures that the jet isomorphism extends to
an isomorphism
$\Phi^\sharp:
\mathcal J^\infty_{\mathrm{lin},\mathcal M}
\longrightarrow
\mathcal J^\infty_{\mathcal M}$.
Let $\widehat\pi_{\mathcal A}$ denote the quotient map
\[
\widehat\pi_{\mathcal A}:
\mathcal J^\infty_{\mathrm{lin},\mathcal M}
\longrightarrow
\mathcal J^\infty_{\mathrm{lin},\mathcal M}/\widehat{\mathcal A}
\simeq
\mathcal O_{\mathcal M}[[\dot z]],
\]
where $\widehat{\mathcal A}$ is the ideal defined in
Section~\ref{sec:ideal-A}.
We define the formal exponential map associated with the higher jet construction by
\begin{equation}\label{eq jet formal exponential}
\exp^\sharp_{\mathrm{jet}}
:=
\widehat\pi_{\mathcal A}
\circ(\Phi^\sharp)^{-1}
\circ\tau^\sharp:
\mathcal O_{\mathcal M}
\longrightarrow
\mathcal O_{\mathcal M}[[\dot z]].
\end{equation}
By Theorem~\ref{thm inverse jets mod A}, for every
$f\in\mathcal O_{\mathcal M}$ we have
\begin{equation}\label{eq jet exponential coordinates}
\exp^\sharp_{\mathrm{jet}}(f)
=
\sum_{n=0}^{\infty}\frac{1}{n!}(\nabla^*)^n(f).
\end{equation}

\begin{proposition}\label{prop:normalized-jet-exponential}
Assume that the Ehresmann connection preserves the zero section.
Then the map
\[
\exp^\sharp_{\mathrm{jet}}
=
\widehat\pi_{\mathcal A}\circ
(\Phi^\sharp)^{-1}\circ\tau^\sharp
\]
is a normalized formal exponential map.
\end{proposition}

\begin{proof}
The canonical Taylor morphism $\tau^\sharp$ is a morphism
of sheaves of superalgebras. By
Corollary~\ref{cor theorem linear}, the morphism $\Phi^\sharp$ is an
isomorphism of sheaves of supercommutative algebras, and
$\widehat\pi_{\mathcal A}$ is the canonical quotient morphism.
Consequently, their composition $\exp^\sharp_{\mathrm{jet}}$
is a morphism of sheaves of supercommutative algebras.
Since each of the morphisms in this composition is locally
determined by the images of the local coordinates, the same holds for $\exp^\sharp_{\mathrm{jet}}$.
Thus,  it suffices to verify the required
conditions on local coordinates $z^a$.

Since the Ehresmann connection preserves the zero section,
its local coefficients satisfy $\Gamma_a^b(z,0)=0$.
The corresponding derivation is
\[
\nabla^*=
\dot z^a\frac{\partial}{\partial z^a}
-\dot z^a\Gamma_a^b(z,\dot z)
\frac{\partial}{\partial\dot z^b}.
\]
Thus $\nabla^*$ increases the fiber degree by at least one.
In particular,
\[
\nabla^*(z^a)=\dot z^a,
\qquad
(\nabla^*)^n(z^a)=O(\dot z^2),
\quad n\ge2.
\]
Using the formula (\ref{eq jet exponential coordinates}),
we conclude that
\[
\exp^\sharp_{\mathrm{jet}}(z^a)
=z^a+\dot z^a+O(\dot z^2).
\]
\end{proof}

\begin{remark}\label{rem:direct-connection-exponential}
The formal exponential map can also be defined directly
from the connection operator, without using the higher
jet construction. Indeed, consider the operator
\[
\nabla^*
=
\dot z^a\frac{\partial}{\partial z^a}
-
\dot z^a\Gamma_a^c(z,\dot z)
\frac{\partial}{\partial\dot z^c}
\]
on $T\mathcal M$, extended formally to
$\mathcal O_{\mathcal M}[[\dot z]]$ by taking the Taylor
expansions of its coefficients along the zero section.
Let $\mathcal J$ be the ideal generated by the fiber
coordinates $\dot z^a$. If the Ehresmann connection
preserves the zero section, then
$\Gamma_a^c(z,0)=0$, and hence
\[
\nabla^*(\mathcal J^k)\subseteq\mathcal J^{k+1},
\qquad k\geq0.
\]
Consequently, the exponential series
\[
e^{\nabla^*}
=
\sum_{n=0}^{\infty}\frac{1}{n!}(\nabla^*)^n
\]
yields a well-defined automorphism of the completed
superalgebra $\mathcal O_{\mathcal M}[[\dot z]]$.
Restricting this automorphism to
$\mathcal O_{\mathcal M}$ gives the morphism
\begin{equation}\label{eq connection formal exponential}
\exp^\sharp_\nabla:
\mathcal O_{\mathcal M}
\longrightarrow
\mathcal O_{\mathcal M}[[\dot z]],
\qquad
\exp^\sharp_\nabla(f)=e^{\nabla^*}f.
\end{equation}
Since
\[
\nabla^*(z^a)=\dot z^a,
\qquad
(\nabla^*)^n(z^a)\in\mathcal J^2
\quad (n\geq2),
\]
we obtain
\[
\exp^\sharp_\nabla(z^a)
=
z^a+\dot z^a+O(\dot z^2).
\]
Thus $\exp^\sharp_\nabla$ is a normalized formal
exponential map.

By \eqref{eq jet exponential coordinates}, this
direct construction coincides with the formal
exponential map associated with higher jets:
\[
\exp^\sharp_\nabla=\exp^\sharp_{\mathrm{jet}}.
\]
\end{remark}

We conclude this section with a simple example, which demonstrates why
the assumption $\Gamma_a^c(z,0)=0$ is essential for the formal constructions above.

\begin{example}\label{example exponential not formal}
Let $\mathcal M=\mathbb C$ with coordinate $x$, and consider the Ehresmann connection
$$
X=\frac{\partial}{\partial x}
-\frac{\partial}{\partial\dot x}.
$$
Thus $\Gamma(x,\dot x)=1$, and the connection does not preserve the zero section. The corresponding operator is
$$
\nabla^*
=
\dot x\frac{\partial}{\partial x}
-\dot x\frac{\partial}{\partial\dot x}.
$$

We first consider the linearization of the finite jet
spaces. Recall that the modified jet coordinates in the
proof of Theorem~\ref{theor Ehresmann jet space}
are defined recursively by
\[
\widetilde x^{(1)}=x^{(1)},\qquad
\widetilde x^{(n)}
=
\partial\bigl(\widetilde x^{(n-1)}\bigr)
-
\nabla^{*(n-1)}
\bigl(\widetilde x^{(n-1)}\bigr),
\quad n\geq2.
\]
Here $\partial$ is the total jet derivative defined in
\eqref{eq total derivative}, and $\nabla^{*(n-1)}$
is the connection operator acting on the
$(n-1)$-st modified jet coordinates.
For the connection $\Gamma(x,\dot x)=1$, we have
\[
\nabla^{*(n-1)}
\bigl(\widetilde x^{(n-1)}\bigr)
=-x^{(1)}.
\]
Thus the recursion reduces to
$\widetilde x^{(n)}
=
\partial\bigl(\widetilde x^{(n-1)}\bigr)+x^{(1)}$.
Starting with $\widetilde x^{(1)}=x^{(1)}$,
we obtain
\[
\begin{aligned}
\widetilde x^{(2)}
&=\partial(x^{(1)})+x^{(1)}
=x^{(2)}+x^{(1)},\\
\widetilde x^{(3)}
&=\partial(x^{(2)}+x^{(1)})+x^{(1)}\\
&=x^{(3)}+x^{(2)}+x^{(1)}.
\end{aligned}
\]
Continuing inductively gives
\[
\widetilde x^{(n)}
=
x^{(n)}+x^{(n-1)}+\cdots+x^{(1)},
\qquad n\geq2.
\]
Consequently, the inverse transformation is
$$
x^{(1)}= \widetilde x^{(1)}, \qquad
x^{(n)}
=
\widetilde x^{(n)}-\widetilde x^{(n-1)},
\quad n\geq2.
$$
Setting $\widetilde x^{(n)}=0$ for $n>1$, we obtain
$$
x^{(1)}=\dot x,\qquad
x^{(2)}=-\dot x,\qquad
x^{(n)}=0,\quad n\geq3.
$$
Thus, applying the inverse linearization and the quotient
by $\mathcal A$ term by term to the Taylor expression gives
$$
x+\sum_{n=1}^{\infty}\frac{x^{(n)}}{n!}
\longmapsto
x+\frac12\dot x.
$$

On the other hand,
$\nabla^*x=\dot x$ and
$\nabla^*\dot x=-\dot x$,
so
$$
(\nabla^*)^n x=(-1)^{n-1}\dot x,
\qquad n\geq1.
$$
The corresponding exponential expression is
$$
e^{\nabla^*}x
=
x+\sum_{n=1}^{\infty}
\frac{(-1)^{n-1}}{n!}\dot x.
$$
This is not a well-defined formal sum in the fiber coordinates, since infinitely many summands contribute to the coefficient of $\dot x$. Nevertheless, the numerical series converges, giving
$$
x+(1-e^{-1})\dot x.
$$
Since
$$
\frac12\neq1-e^{-1},
$$
the termwise jet construction and the analytically summed
exponential expression give different results.

This example illustrates why the zero-section-preserving
assumption is needed: without it, neither the extension of
the jet linearization to the algebra of formal series nor
the formal exponential of $\nabla^*$ is guaranteed to exist
in this algebra.
\end{example}

\begin{remark}\label{rem:inverse-linearization-not-injective}
The inverse coordinate transformation
\[
x^{(1)}=\widetilde x^{(1)},\qquad
x^{(n)}
=
\widetilde x^{(n)}-\widetilde x^{(n-1)},
\quad n\geq2,
\]
extends term by term to infinite jet series. However,
this extension is not injective.
Indeed, consider the nonzero formal jet series
\[
F=\sum_{n=1}^{\infty}x^{(n)}.
\]
Applying out map, we obtain
\[
\begin{aligned}
F
&\longmapsto
\widetilde x^{(1)}
+\sum_{n=2}^{\infty}
\bigl(\widetilde x^{(n)}
-\widetilde x^{(n-1)}\bigr)=0,
\end{aligned}
\]
since the coefficient of every
$\widetilde x^{(n)}$ vanishes.

Thus the inverse linearization, although invertible
at every finite jet order, is not injective on the
space of infinite jet series. In particular, the
finite-order isomorphisms do not automatically extend
to an isomorphism of the corresponding completed
jet spaces.
\end{remark}

\section{A Koszul-type splitting theorem for supermanifolds}

\subsection{A splitting theorem for formal exponential maps}

Let $\mathcal M$ be a supermanifold and let
\[
F^\sharp:
\mathcal O_{\mathcal M}
\longrightarrow
\mathcal O_{\mathcal M}[[\dot z]]
\]
be a formal exponential map. We show that $F^\sharp$
determines an isomorphism between $\mathcal M$ and
its associated graded supermanifold.
Choose local even and odd coordinates
$(x^i,\xi^\alpha)$ on $\mathcal M$, with corresponding
fiber coordinates $(\dot x^i,\dot\xi^\alpha)$ on
$T\mathcal M$. Consider the ideal
\[
\mathcal I
=
\langle \xi^\alpha,\dot x^i\rangle
\subset
\mathcal O_{\mathcal M}[[\dot z]].
\]
This ideal is globally defined. Indeed, under a
change of coordinates
\[
x^{i'}=x^{i'}(x,\xi),
\qquad
\xi^{\alpha'}=\xi^{\alpha'}(x,\xi),
\]
the odd coordinates $\xi^{\alpha'}$ belong to the
ideal generated by the $\xi^\alpha$, while the
transformation law for the even tangent coordinates is
\[
\dot x^{i'}
=
\dot z^a
\frac{\partial x^{i'}}{\partial z^a}.
\]
Modulo $\mathcal I$, the terms containing
$\dot x^i$ vanish, and the coefficients of the
remaining terms contain odd base coordinates.
Thus $\dot x^{i'}\in\mathcal I$, proving that
$\mathcal I$ is independent of the choice of
local coordinates. A related construction was considered in \cite{Vi}.

We have
\[
\mathcal O_{\mathcal M}[[\dot z]]/\mathcal I
\simeq
\mathcal O_{\operatorname{gr}\mathcal M}.
\]
Locally, the quotient is obtained by setting
$\xi^\alpha=0$ and $\dot x^i=0$, leaving the
even coordinates $x^i$ and the odd coordinates
$\dot\xi^\alpha$.
Denote the quotient morphism by
\[
i^\sharp:
\mathcal O_{\mathcal M}[[\dot z]]
\longrightarrow
\mathcal O_{\operatorname{gr}\mathcal M}.
\]
Composing it with $F^\sharp$, we obtain
\[
s_F^\sharp
:=
i^\sharp\circ F^\sharp:
\mathcal O_{\mathcal M}
\longrightarrow
\mathcal O_{\operatorname{gr}\mathcal M}.
\]

\begin{theorem}\label{thm formal exponential splitting}
Every formal exponential map $F^\sharp$ induces a
splitting of $\mathcal M$. More precisely, the morphism
\[
s_F^\sharp=i^\sharp\circ F^\sharp:
\mathcal O_{\mathcal M}
\xrightarrow{\sim}
\mathcal O_{\operatorname{gr}\mathcal M}
\]
is an isomorphism of sheaves of supercommutative
algebras.
\end{theorem}

\begin{proof}
By the definition of a formal exponential map in
Section~\ref{sec:formal-exponential-maps}, in local
coordinates $z^a=(x^i,\xi^\alpha)$ we have
\[
F^\sharp(z^a)
=
z^a+A^a_{b}(z)\dot z^b+O(\dot z^2),
\]
where $A(z)$ is an invertible even supermatrix. Passing to the quotient by $\mathcal I$, we set
$\xi^\alpha=0$ and $\dot x^i=0$. Since $F^\sharp$
preserves parity, the resulting expressions have
the form
\[
s_F^\sharp(x^i)
=
x^i+O(\dot\xi^2),
\qquad
s_F^\sharp(\xi^\alpha)
=
B^\alpha_{\beta}(x)\dot\xi^\beta
+O(\dot\xi^3),
\]
where
\[
B^\alpha_{\beta}(x)
=
A^\alpha_{\beta}(x,0).
\]
Here $O(\dot\xi^k)$ denotes terms of degree at least
$k$ in the odd coordinates $\dot\xi$.

Since $A(z)$ is an invertible even supermatrix
and $B(x)$ is the odd--odd block of its reduction
to the reduced manifold, $B(x)$ is invertible.  It follows that $s_F^\sharp$ induces the identity
on the reduced manifold and an isomorphism on the
first-order odd neighborhoods.
The local expressions above show that $s_F^\sharp$
is a local isomorphism, since $B(x)$ is invertible.
Moreover, $s_F^\sharp$ induces the identity on the
underlying reduced manifold. 
Therefore $s_F^\sharp$ is locally an isomorphism
of sheaves of supercommutative algebras.
Since $s_F^\sharp$ is globally defined and induces
the identity on the underlying reduced manifold,
the local inverses agree on overlaps and define
a global inverse. Hence
$\mathcal O_{\mathcal M}
\simeq
\mathcal O_{\operatorname{gr}\mathcal M}$,
which gives a splitting of $\mathcal M$.
\end{proof}

If $F^\sharp$ is normalized, then
$A^a_{\ b}(z)=\delta^a_b$. Consequently,
$B^\alpha_{\ \beta}(x)=\delta^\alpha_\beta$, and
the induced splitting has the local expression
\[
s_F^\sharp(x^i)
=
x^i+O(\dot\xi^2),
\qquad
s_F^\sharp(\xi^\alpha)
=
\dot\xi^\alpha+O(\dot\xi^3).
\]

\begin{corollary}\label{cor splitting from connection}
Every zero-section-preserving Ehresmann connection
on $\mathcal M$ induces a splitting of $\mathcal M$.
\end{corollary}

\begin{proof}
By
Proposition~\ref{prop:normalized-jet-exponential}
and Remark~\ref{rem:direct-connection-exponential},
the formal exponential maps
$\exp^\sharp_\nabla$ and
$\exp^\sharp_{\mathrm{jet}}$ coincide and induce
the same splitting, described by the commutative
diagram
\[
\xymatrix{
& \mathcal O_{\mathcal M}[[\dot z]]
  \ar[rd]^{i^\sharp} & \\
\mathcal O_{\mathcal M}
  \ar[ur]^{\exp^\sharp_\nabla
           =\exp^\sharp_{\mathrm{jet}}}
  \ar[rr]_{\simeq}^{s^\sharp_\nabla
                   =s^\sharp_{\mathrm{jet}}}
&&
\mathcal O_{\operatorname{gr}\mathcal M}.
}
\]
\end{proof}

\subsection{The role of the zero-section-preserving condition}

By Corollary~\ref{cor splitting from connection}, an
Ehresmann connection preserving the zero section
determines two coinciding formal exponential maps,
which induce the same splitting of $\mathcal M$.  Without the zero-section-preserving assumption, the two
constructions need not be well defined in the algebra of
formal power series, since infinitely many terms may
contribute to the same coefficient. Nevertheless, in some
cases these coefficients are given by convergent numerical
series and can be evaluated by analytic summation.

The resulting expressions may then exhibit four different
behaviors: both yield splittings, only the jet construction
yields a splitting, only the direct construction yields a
splitting, or neither yields a splitting.
We illustrate the first three possibilities in a single
family of examples. Throughout these examples, the infinite
sums are understood in the sense of analytic summation of
their coefficients.

\begin{example}\label{example jet not splitting}
Let $\mathcal M=\mathbb C^{0|1}$ with odd coordinate
$\xi$, and consider the Ehresmann connection
$$
X=\frac{\partial}{\partial\xi}
-\lambda\frac{\partial}{\partial\dot\xi},
\qquad \lambda\in\mathbb C\setminus\{0\}.
$$
The corresponding connection operator is
$$
\nabla^*
=
\dot\xi\frac{\partial}{\partial\xi}
-\lambda\dot\xi\frac{\partial}{\partial\dot\xi}.
$$

\medskip
\noindent
\textbf{The jet construction.}
Recall that the recursive construction in the proof of
Theorem~\ref{theor Ehresmann jet space} gives
\[
\widetilde\xi^{(1)}=\xi^{(1)},\qquad
\widetilde\xi^{(k+1)}
=
\partial\widetilde\xi^{(k)}
-\nabla^{*(k)}\widetilde\xi^{(k)}.
\]
For the connection under consideration, we have
$\nabla^{*(k)}\widetilde\xi^{(k)}
=-\lambda\xi^{(1)}$.
Hence the recursion reduces to
$\widetilde\xi^{(k+1)}
=
\partial\widetilde\xi^{(k)}+\lambda\xi^{(1)}$.
In particular,
\[
\begin{aligned}
\widetilde\xi^{(2)}
&=\xi^{(2)}+\lambda\xi^{(1)},\qquad
\widetilde\xi^{(3)}
&=\xi^{(3)}+\lambda\xi^{(2)}+\lambda\xi^{(1)}.
\end{aligned}
\]
Proceeding inductively, we obtain
\[
\widetilde\xi^{(k)}
=
\xi^{(k)}
+\lambda\sum_{j=1}^{k-1}\xi^{(j)},
\qquad k\geq2.
\]
Passing to the quotient by
$\mathcal A=\langle\widetilde\xi^{(k)};\ k>1\rangle$,
we get
\[
\xi^{(k)}
=-\lambda\sum_{j=1}^{k-1}\xi^{(j)},
\qquad k\geq2.
\]
In particular, $\xi^{(2)}=-\lambda\xi^{(1)}$.
Subtracting two consecutive relations gives
\[
\xi^{(k+1)}=(1-\lambda)\xi^{(k)},
\qquad k\geq2.
\]
Consequently, in the quotient we have
\[
\xi^{(1)}=\dot\xi,
\qquad
\xi^{(k)}
=-\lambda(1-\lambda)^{k-2}\dot\xi,
\quad k\geq2.
\]

Substituting these relations
into $\sum\limits_{k=1}^{\infty}\frac{\xi^{(k)}}{k!}$ 
we obtain
\[
\begin{aligned}
s_{\mathrm{jet}}^\sharp(\xi)
&=
\dot\xi
-\frac{\lambda}{2!}\dot\xi
-\frac{\lambda(1-\lambda)}{3!}\dot\xi
-\frac{\lambda(1-\lambda)^2}{4!}\dot\xi
-\cdots\\
&=
\left(
1-\lambda\sum_{k=2}^{\infty}
\frac{(1-\lambda)^{k-2}}{k!}
\right)\dot\xi=C(\lambda)\dot\xi,
\end{aligned}
\]
where
\[
\begin{aligned}
C(\lambda)
&=
1-\lambda
\frac{e^{1-\lambda}-1-(1-\lambda)}
{(1-\lambda)^2}\\[4pt]
&=
\frac{1-2\lambda+\lambda^2-\lambda e^{1-\lambda}
+\lambda+\lambda-\lambda^2}
{(1-\lambda)^2}=
\frac{1-\lambda e^{1-\lambda}}
{(1-\lambda)^2}.
\end{aligned}
\]
The expression at $\lambda=1$ is understood
by continuity.

\medskip
\noindent
\textbf{The direct construction.}
Since
$\nabla^*\xi=\dot\xi$,
$\nabla^*\dot\xi=-\lambda\dot\xi$,
we have
$$
(\nabla^*)^n\xi
=
(-\lambda)^{n-1}\dot\xi,
\qquad n\geq1.
$$
Further,
\[
\begin{aligned}
s_\nabla^\sharp(\xi)
&=
\sum_{n=1}^{\infty}
\frac{1}{n!}(\nabla^*)^n\xi=
\dot\xi-\frac{\lambda}{2!}\dot\xi
+\frac{\lambda^2}{3!}\dot\xi
-\frac{\lambda^3}{4!}\dot\xi+\cdots\\
&=
\left(
\sum_{n=1}^{\infty}
\frac{(-\lambda)^{n-1}}{n!}
\right)\dot\xi=
D(\lambda)\dot\xi,\quad D(\lambda)=\frac{1-e^{-\lambda}}{\lambda},
\end{aligned}
\]
with $D(0)=1$ by continuity.

Both constructions therefore define morphisms
from $\mathcal O_{\mathcal M}$ to
$\mathcal O_{\operatorname{gr}\mathcal M}$.
Since $\mathcal M$ is purely odd of dimension
$0|1$, these morphisms are isomorphisms precisely
when their respective coefficients are nonzero.

\medskip
\noindent
\textbf{Case 1. Both constructions give splittings.}
Take $\lambda=1$. Then
$$
C(1)=\frac12,
\qquad
D(1)=1-e^{-1}.
$$
Both coefficients are nonzero, so both
constructions give splittings. However,
$s_{\mathrm{jet}}^\sharp
\neq s_\nabla^\sharp$.
Thus the two splittings need not coincide.

\medskip
\noindent
\textbf{Case 2. Only the jet construction gives
a splitting.}
Take $\lambda=2\pi i$. Then
$$D(2\pi i)=\frac{1-e^{-2\pi i}}{2\pi i}=0.
$$
Since $e^{1-2\pi i}=e$, we obtain
\[
\begin{aligned}
C(2\pi i)
&=\frac{1-2\pi i e}{(1-2\pi i)^2}
\neq0.
\end{aligned}
\]
Consequently, $s_{\mathrm{jet}}^\sharp$
is an isomorphism, while $s_\nabla^\sharp$
is not.

\medskip
\noindent
\textbf{Case 3. Only the direct construction gives
a splitting.}
At any zero $\lambda\neq1$ of $C$, we have
$\lambda e^{-\lambda}=e^{-1}$, which has infinitely many complex solutions,
given by $\lambda=-W_k(-e^{-1})$, where $W_k$
denotes a branch of the Lambert $W$-function. (The Lambert $W$-function is the multivalued inverse
of the function $w\mapsto we^w$. Its branches
$W_k$, $k\in\mathbb Z$, satisfy
$W_k(z)e^{W_k(z)}=z.$) 
Hence
\[
D(\lambda)
=\frac{1-e^{-\lambda}}{\lambda}
=\frac{e\lambda-1}{e\lambda^2}.
\]
This is nonzero, since $\lambda=e^{-1}$ does not
satisfy $\lambda e^{-\lambda}=e^{-1}$.
Thus $s_\nabla^\sharp$ defines a splitting,
whereas $s_{\mathrm{jet}}^\sharp$ does not.
\end{example}

\begin{example}\label{example neither splitting}
The fourth possibility can be realized on
$\mathcal M=\mathbb C^{0|2}$.
Consider the Ehresmann connection
$$
X_\alpha
=
\frac{\partial}{\partial\xi^\alpha}
-\lambda_\alpha
\frac{\partial}{\partial\dot\xi^\alpha},
\qquad \alpha=1,2.
$$
The two coordinate directions are independent,
so the preceding calculations give
$$
s_{\mathrm{jet}}^\sharp(\xi^\alpha)
=
C(\lambda_\alpha)\dot\xi^\alpha,
\qquad
s_\nabla^\sharp(\xi^\alpha)
=
D(\lambda_\alpha)\dot\xi^\alpha.
$$
Choose
$$
\lambda_1=2\pi i,
\qquad
\lambda_2=\lambda_0,
$$
where $C(\lambda_0)=0$, as in
Example~\ref{example jet not splitting}.

Since
$$
D(\lambda_1)=0,
\qquad
C(\lambda_2)=0,
$$
both matrices of linear terms are singular.
Hence neither construction defines a splitting.
\end{example}

\bigskip
\bigskip

\noindent
L.~C.: Departamento de Matem{\'a}tica, Instituto de Ci{\^e}ncias Exatas,
Universidade Federal de Minas Gerais, BRAZIL.

\medskip

\noindent Email: {\tt luizfelipe.ac@hotmail.com}
\bigskip

\noindent A.~Sh.:   Department of Quantum Field Theory,    Physics Faculty, Tomsk State University, RUSSIA.

\medskip

\noindent Email: {\tt sharapov@phys.tsu.ru}

\bigskip

\noindent
E.~V.: Departamento de Matem{\'a}tica, Instituto de Ci{\^e}ncias Exatas,
Universidade Federal de Minas Gerais, BRAZIL.

\medskip

\noindent Email: {\tt VishnyakovaE\symbol{64}googlemail.com}


\begin{thebibliography}{99}

\bibitem{Ber} {\it F. A. Berezin.} Introduction to algebra and analysis with anticommuting variables.
V. P. Palamodov, ed., Moscow State University Press, Moscow, 1983. Expanded transl.
into English: Introduction to superanalysis. A. A. Kirillov, ed., D. Reidel, Dordrecht, 1987.



\bibitem{B}
A.~Bruce,
\textit{On curves and jets on supermanifolds},
Arch. Math. (Brno), 2014.

\bibitem{CF}
A.~S.~Cattaneo and G.~Felder,
\textit{On the globalization of Kontsevich's star product
and the perturbative Poisson sigma model},
Progr. Theoret. Phys. Suppl. \textbf{144}
(2001), 38--53.

\bibitem{Correa2026-gm}
M.~Corrêa and S.~Noja,
\textit{Formal moduli and the splitting theory of supermanifolds},
arXiv:2605.03166 [math.AG] (2026).

\bibitem{LS} Hsuan-Yi Liao and Mathieu Stiénon. ``Formal exponential map for graded manifolds." International Mathematics Research Notices 2019.3 (2019): 700-730; arXiv:1508.02780 [math.DG]. 

\bibitem{Emmrich1993TheDG}
C. Emmrich and A. Weinstein. ``The differential geometry of Fedosov’s quantization." Lie Theory and Geometry: In Honor of Bertram Kostant. Boston, MA: Birkhäuser Boston, 1994. 217-239; 	arXiv:hep-th/9311094. 




\bibitem{Leites}
D.~A.~Leites,
\textit{Introduction to the theory of supermanifolds},
Russian Math. Surveys \textbf{35} (1980), no.~1,
1--64.

\bibitem{Kos}
J.-L.~Koszul,
\textit{Connections and splittings of supermanifolds},
Differential Geom. Appl. \textbf{4} (1994), 151--161.

\bibitem{Vi}
E.~Vishnyakova,
\textit{Graded covering of a supermanifold},
arXiv:2212.09558.






\end{thebibliography}
\end{document}